# The *q*-Levenberg-Marquardt method for unconstrained nonlinear optimization

Danijela D. Protić
Center for Applied Mathematics and Electronic
Belgrade 11090
adanijela@ptt.rs

Miomir S. Stanković
Mathematical Institute of SASA
Belgrade 11000
miomir.stankovic@gmail.com

## ABSTRACT

A *q*-Levenberg-Marquardt method is an iterative procedure that blends a *q*-steepest descent and *q*-Gauss-Newton methods. When the current solution is far from the correct one the algorithm acts as the *q*-steepest descent method. Otherwise the algorithm acts as the *q*-Gauss-Newton method. A damping parameter is used to interpolate between these two methods. The *q*-parameter is used to escape from local minima and to speed up the search process near the optimal solution.

## 1 Introduction

Nonlinear optimization algorithms minimize the objective function by fitting a parametrized model to a set of data points [1], [2]. The steepest-descent (SD) method minimizes the objective function by an iterative procedure based on the step length and search direction determined by the negative of the gradient [3]. The SD uses the first-order derivatives of the objective function and quickly approaches the solution from a distance. However, its convergence becomes slow near the solution [4]. The Gauss-Newton (GN) method presumes that the error function is approximately quadratic near the optimal solution. The method simplifies the process of calculating the Hessian matrix directly with the Jacobian matrix introduced and by minimizing the sum of the squares of the errors (SSE) [5]. The Levenberg-Marquardt (LM) method acts like the SD method when the parameters are far from their optimal value and acts like the GN method when parameters are close to their optimal value [4].

Over the years, various methods for solving nonlinear optimization problems based on the *q*-calculus have been proposed. The *q*-steepest-descent (*q*-SD) algorithms are presented in [3], [6-8]. The algorithms are iterative and solve nonlinear optimization problems with various strategies for generating the *q*-parameters and computing the step lengths. The *q*-Gauss-Newton (*q*-GN) algorithm converges only when the iteration process starts near to the optimal solution and speeds up the convergence near the global minimum [9]. The algorithm does not require evaluations of *q*-second-order derivatives.

In this paper we present a new LM model based on *q*-calculus. The *q*-Levenberg-Marquardt (*q*-LM) algorithm is an iterative procedure which uses damping factor $\lambda$ to interpolate between *q*-SD and *q*-GN algorithms. The *q*-LM algorithm minimizes the distance between current and previous iterate as well as the SSE of the corresponding values of the objective function. If $\lambda \to 0$ the algorithm acts as the *q*-GN method, which estimates of the *q*-first-order derivatives of the objective function and speeds up the convergence near the stationary point. If $\lambda \to \infty$ the algorithm tends to the *q*-SD algorithm which quickly approaches the solution from the distance.

## 2 Levenberg-Marquardt method

Suppose that $f: \mathbb{R}^n \to \mathbb{R}$ is a real valued continuous function taking as input vector $\mathbf{x} = [x_1 \ x_2 \cdots x_n]^{\mathrm{T}}$. The goal of the LM method is to minimize Euclidean distance given in the form $\|f(x)\|^2$. The gradient $\nabla f(\mathbf{x})$ and the Hessian $\mathbf{H}(f(\mathbf{x}))$ are vectors of partial derivatives and matrix of second partial derivatives of $f(\mathbf{x})$, respectively, as follows

$$\nabla f(\mathbf{x}) = Df(\mathbf{x}) = \left[\frac{\partial f(\mathbf{x})}{\partial x_1} \frac{\partial f(\mathbf{x})}{\partial x_2} \cdots \frac{\partial f(\mathbf{x})}{\partial x_n}\right], \quad \mathbf{H}\big(f(\mathbf{x})\big) = \nabla^2 f(\mathbf{x}) = D^2 f(\mathbf{x}) = \begin{bmatrix} \frac{\partial^2 f(\mathbf{x})}{\partial x_1^2} & \frac{\partial^2 f(\mathbf{x})}{\partial x_1 \partial x_2} & \cdots & \frac{\partial^2 f(\mathbf{x})}{\partial x_1 \partial x_n} \\ \frac{\partial^2 f(\mathbf{x})}{\partial x_2 \partial x_1} & \frac{\partial^2 f(\mathbf{x})}{\partial x_2^2} & \cdots & \frac{\partial^2 f(\mathbf{x})}{\partial x_2 \partial x_n} \\ \vdots & & \ddots & \vdots \\ \frac{\partial^2 f(\mathbf{x})}{\partial x_n \partial x_1} & \frac{\partial^2 f(\mathbf{x})}{\partial x_n \partial x_2} & \cdots & \frac{\partial^2 f(\mathbf{x})}{\partial x_n^2} \end{bmatrix}_{n \times n} \quad (1)$$

In high dimensional space for any $g_i = \frac{\partial f_j(\mathbf{x})}{\partial x_i}$ ($i = 1, \dots, n$; $j = 1, \dots m$), the gradient of partial derivatives with respect to each dimension is given with $g_j(\mathbf{x}) = \nabla f_j(\mathbf{x}) = [g_{1,j}(\mathbf{x})\, g_{2,j}(\mathbf{x}) \cdots g_{n,j}(\mathbf{x})]^\mathrm{T}$. The Jacobian matrix of *f*, $\mathbf{J}(f(x))$, contains all its first-order partial derivatives [10]

$$\mathbf{J}(f(x)) = \begin{bmatrix} \frac{\partial f_1(\mathbf{x})}{\partial x_1} & \frac{\partial f_1(\mathbf{x})}{\partial x_2} & \cdots & \frac{\partial f_1(\mathbf{x})}{\partial x_n} \\ \frac{\partial f_2(\mathbf{x})}{\partial x_1} & \frac{\partial f_2(\mathbf{x})}{\partial x_2} & \cdots & \frac{\partial f_2(\mathbf{x})}{\partial x_n} \\ \vdots & & \ddots & \vdots \\ \frac{\partial f_m(\mathbf{x})}{\partial x_1} & \frac{\partial f_m(\mathbf{x})}{\partial x_2} & \cdots & \frac{\partial f_m(\mathbf{x})}{\partial x_n} \end{bmatrix}_{n \times m} \quad (2)$$

The Hessian matrix can be determined as: $\mathbf{H}(f(\mathbf{x})) = \mathbf{J}(\nabla f(\mathbf{x}))^\mathrm{T}$.

Consider the nonlinear function *f*: R→R. In solving nonlinear least squares (NLS) problems the idea is to form an affine transformation of *f* near the current point. The truncated first-order Taylor approximation of *f* at any point *p* is given with the following formul

$$\hat{f}(x; p) = f(p) + Df(p)(x - p) \quad (3)$$

where $Df(p)(x - p)$ represents the Jacobian matrix of the partial derivatives of *f*. If $\hat{f}(x; p) = f(x)$, provided *x* is near *p*, $\hat{f}(x; p)$ can be minimized using linear least squares (LS) to find exact solution [11].

The LM algorithm is an iterative procedure with current iterate approximate the problem in such extent that $\hat{f}(x; p) = f(x)$ but also $x \approx p$. Let $x^{(1)}, x^{(2)}, \dots x^{(l)}$ be the iterates. At iteration *k* the affine transformation of *f* at point *k* is given with

$$\hat{f}(x; x^{(k)}) = f(x^{(k)}) + Df(x^{(k)})(x - x^{(k)}) \quad (4)$$

The idea is to choose $x^{(k+1)}$ as a minimizer of

$$\|\hat{f}(x; x^{(k)})\|^2 + \lambda^{(k)} \|x - x^{(k)}\|^2, \quad \lambda^{(k)} > 0 \quad (5)$$

The goal is to minimize first and second term in (5). If $\|\hat{f}(x; x^{(k)})\|^2$ is minimized then the approximation can be considered as $\hat{f} \approx f$. Second objective determines how far $x^{(k)}$ is from $x$. The LM algorithm describes the trade between these two objectives. Parameter $\lambda^{(k)}$ varies with the step size and gives the level of distrust. Suppose $x^{(k+1)}$ is a solution of the LS problem

$$\text{minimize } \left\{ \|f(x^{(k)}) + Df(x^{(k)})(x - x^{(k)})\|^2 + \lambda^{(k)} \|x - x^{(k)}\|^2 \right\}, \quad \lambda^{(k)} > 0 \quad (6)$$

then

$$x^{(k+1)} = x^{(k)} - \left(Df(x^{(k)})^\mathrm{T} Df(x^{(k)}) + \lambda^{(k)} \mathrm{I}\right)^{-1} Df(x^{(k)})^\mathrm{T} f(x^{(k)}), \quad \lambda^{(k)} > 0 \quad (7)$$

the term $\left(Df(x^{(k)})^\mathrm{T} Df(x^{(k)}) + \lambda^{(k)} \mathrm{I}\right)^{-1}$ represents a general solution of bi-criterion and always exists since $\lambda^{(k)} > 0$. The point $x^{(k)}$ is a stationary point if $x^{(k+1)} \approx x^{(k)}$. Only if $Df(x^{(k)})^\mathrm{T} f(x^{(k)}) = 0$ the optimal condition holds. During the iterative procedure $\lambda^{(k)}$ is adjusted. If $\lambda^{(k)}$ is too big $x^{(k+1)}$ is too close to $x^{(k)}$, and progress is slow. Otherwise $x^{(k+1)}$ is far from $x^{(k)}$ and approximation is poor. The update mechanism is given as follows:

1. if $\|f(x^{(k+1)})\|^2 < \|f(x^{(k)})\|^2$, $f(x^{(k+1)})$ is better than the current objective accept new $x^{(k+1)}$ and reduce $\lambda^{(k)}$,
2. otherwise increase $\lambda^{(k)}$ and do not update $x^{(k)}$; $x^{(k+1)} = x^{(k)}$.

Consider the iterate $x^{(k+1)}$ as the solution of the gradient of the objective function (the SD algorithm)

$$x^{(k+1)} = x^{(k)} - \lambda^{(k)} Df(x^{(k)}), \quad \lambda^{(k)} > 0 \qquad (8)$$

and the iterate $x^{(k+1)}$ as the solution of the GN algorithm near the optimal value

$$x^{(k+1)} = x^{(k)} - \left(D^2 f(x^{(k)})\right)^{-1} Df(x^{(k)})^{\mathrm{T}} \qquad (9)$$

The LM algorithm starts from initial point $x^{(0)}$ and initial factor $\lambda^{(0)}$. Afterwards the iterate $x^{(k+1)}$ can be determined with the formula

$$x^{(k+1)} = x^{(k)} - \left(\mathbf{H} + \lambda^{(k)} \mathbf{I}\right)^{-1} \mathbf{J}^{\mathrm{T}} f(x^{(k)}) = x^{(k)} - \left(\mathbf{J}^{\mathrm{T}} \mathbf{J} + \lambda^{(k)} \mathbf{I}\right)^{-1} \mathbf{J}^{\mathrm{T}} f(x^{(k)}) \qquad (10)$$

where $\mathbf{I}$ represents the identity matrix and $\mathbf{H} = \mathbf{J}^{\mathrm{T}} \mathbf{J}$. If $\lambda^{(k)} \to \infty$, then $\mathbf{H} + \lambda^{(k)} \mathbf{I} \approx \mathbf{I}$, and the LM algorithm acts as the SD algorithm. If $\lambda^{(k)} \to 0$ the LM algorithm acts as the GN algorithm since $x^{(k)}$ is close to the optimum value.

## 3 q-Levenberg-Marquardt method

Consider function $f(x): \mathrm{R} \to \mathrm{R}, x \in \mathrm{R}$. The q-derivative of f is given with the formula [12], [13]:

$$D_q f(x) = \frac{f(x) - f(qx)}{(1-q)x}, \quad x \neq 0, 0 \leq q < 1 \qquad (11)$$

As $q \to 1$ the q-derivative tends to classical derivative [8], [12]. The characteristics of the q-derivative are described in details in [3], [6], [12], [13], [15]. Suppose the partial derivatives of $f_j$: $\mathrm{R}^n \to \mathrm{R}^m$, $m \geq n$, ($i=1,\ldots,n$; $j=1,\ldots,m$) exist. In that case the q-partial derivative of $f_j$, with respect to $x_i$ is

$$D_{q,x_i} f_j(x) = \begin{cases} \frac{f_j(x) - \varepsilon_{q,i} f_j(x)}{(1-q)x_i} & x_i \neq 0 \\ \frac{\partial f_j}{\partial x_i} & x_i = 0 \end{cases} \qquad (12)$$

For parameter vector $\mathbf{q} = [q_1\ q_2\ \cdots\ q_n]^{\mathrm{T}}$ and function $f(\mathbf{x})$ of n variables $\mathbf{x} = [x_1\ x_2\ \cdots\ x_n]^{\mathrm{T}}$, the q-gradient vector can be determined as $\nabla_q f(\mathbf{x}) = \left[D_{q_1,x_1} f(\mathbf{x})\ \cdots\ D_{q_i,x_i} f(\mathbf{x})\ \cdots\ D_{q_n,x_n} f(\mathbf{x})\right]$.

In [12] authors derived the q-analog of the truncated Taylor expansion for any polynomial f(x) of degree n and any number c given with the following formula:

$$f(x) = \sum_{j=0}^{n} \left(D_q^j f\right)(c) \frac{(x-c)_q^j}{j!} \qquad (13)$$

where $D_q^j f$ represents j-th order q-derivative of the objective function f. The q-Jacobian matrix of f is given with

$$D_{q,\mathbf{x}} f(\mathbf{x}) = \left[D_{q,x_j} f_i(\mathbf{x})\right]_{m \times n}, \quad i = 1, \ldots, n;\ j = 1, \ldots, m \qquad (14)$$

At the position i,j the q-Jacobian elements are given as follows

$$\mathbf{J}_q(\mathbf{x})_{i,j} = \frac{f_i(x_1,\cdots,x_{j-1},x_j,x_{j+1},\cdots x_n) - f_i(x_1,\cdots,x_{j-1},qx_j,x_{j+1},\cdots x_n)}{(1-q)x_j} \qquad (15)$$

Consider real-valued function $f(\mathbf{x}): \mathrm{R} \to \mathrm{R}^n$ of $\mathbf{x} = [x_1\ x_2\ \cdots\ x_n]^{\mathrm{T}}$. The basic idea of the q-LM algorithm is to minimize the expression:

$$\left\| f(\mathbf{x}^{(k)}) + D_q f(\mathbf{x}^{(k)})(\mathbf{x} - \mathbf{x}^{(k)}) \right\|^2 + \lambda^{(k)} \left\| \mathbf{x} - \mathbf{x}^{(k)} \right\|^2, \quad \lambda^{(k)} > 0 \qquad (16)$$

At each algorithm iteration, the iterate $\mathbf{x}^{(k+1)}$ is determined according to the formula:

$$\mathbf{x}^{(k+1)} = \mathbf{x}^{(k)} - D_q f(\mathbf{x}^{(k)})^{\dagger} f(\mathbf{x}^{(k)}) \qquad (17)$$

where $D_q f(\mathbf{x}^{(k)})^{\dagger} = \left(D_q f(\mathbf{x}^{(k)})^{\mathrm{T}} D_q f(\mathbf{x}^{(k)}) + \lambda^{(k)} \mathbf{I}\right)^{-1} D_q f(\mathbf{x}^{(k)})^{\mathrm{T}}$, $\lambda^{(k)} > 0$ represents the Moore-Penrose matrix inverse [16], [17]. In accordance with (1), the iterate $\mathbf{x}^{(k+1)}$ can be determined as

$$\mathbf{x}^{(k+1)} = \mathbf{x}^{(k)} - \left(\mathbf{H}_q + \lambda^{(k)} \mathbf{I}\right)^{-1} \mathbf{J}_q^{\mathrm{T}} f(\mathbf{x}^{(k)}), \quad \lambda^{(k)} > 0 \qquad (18)$$

The update mechanism of the $q$-LM algorithm given in the text above terminates only if $D_q f(\mathbf{x}^{(k)})^\mathrm{T} f(\mathbf{x}^{(k)}) = 0$ when the optimal condition holds, the iterate $\mathbf{x}^{(k+1)} \approx \mathbf{x}^{(k)}$, i.e. when the iterate $\mathbf{x}^{(k)}$ becomes the stationary point.

When SSE→0, $\mathbf{x}^{(k+1)} \approx \mathbf{x}^{(k)}$, the $q$-Hessian matrix $\mathbf{H}_q$ can be approximated with $\mathbf{J}_q^\mathrm{T}\mathbf{J}_q$, and $\mathbf{x}^{(k+1)}$ may be determined as follows:

$$\mathbf{x}^{(k+1)} = \mathbf{x}^{(k)} - \left(\mathbf{J}_q^\mathrm{T}\mathbf{J}_q + \lambda^{(k)}\mathbf{I}\right)^{-1}\mathbf{J}_q^\mathrm{T} f(\mathbf{x}^{(k)}) \tag{19}$$

Given update algorithm does not require calculation of second-order $q$-derivatives. The $q$-Jacobian matrix is introduced to simplify the process of calculating the $q$-Hessian matrix directly.

The $q$-LM algorithm can be expressed in the pseudo code as follows:

1. for given function $f$ of $n$ variables $\mathbf{x}$ analytically determine the $q$-Jacobian $\mathbf{J}q$;
2. set:
    {$q$;
    initial guess $\mathbf{x}(0)$;
    initial guess $\lambda(0)$;
    counter $k = 1$;
    *mf* - multiplication factor, *mf*>1;
    *df* - division factor, 0< *df*<1;
    stop iteration $stop\_iter$ to small positive value;
    maximum number of iterations $max\_no\_iter$;};
3. while(1)
    evaluate $f(k)$;
    evaluate $\mathbf{J}q(k)$;
    determine $\mathbf{x}(k+1) = \mathbf{x}(k) - \left(\mathbf{J}_q^\mathrm{T}(k)\mathbf{J}_q(k) + \lambda(k)\mathbf{I}\right)^{-1}\mathbf{J}_q^\mathrm{T}(k)f(\mathbf{x}(k))$;
    3.1. if $\|f(\mathbf{x}(k+1))\|^2 < \|f(\mathbf{x}(k))\|^2$
        accept new $\mathbf{x}(k+1)$;
        reduce $\lambda(k+1) = df \cdot \lambda(k)$;
    3.2. otherwise
        increase $\lambda(k)$; $\lambda(k+1) = mf \cdot \lambda(k)$;
        set $\mathbf{x}(k+1) = \mathbf{x}(k)$;
    if $\|f(\mathbf{x}(k))\|^2 \leq stop\_iter$ OR $k = max\_no\_iter$
        break
    otherwise $k = k + 1$;
    end
end

Since the $q$-LM algorithm requires calculation of the $q$-Jacobian matrix at each iteration step, the formula for the $q$-Jacobian has to be analytically determined before the $q$-LM algorithm starts. Considering processing power, calculating the $q$-Jacobian matrix is less demanding then calculating the $q$-Hessian matrix, which speeds up the algorithm. When $\lambda(k) \to 0$, $\mathbf{x}^{(k)}$ is close to the optimum value and the $q$-LM algorithm acts as the $q$-GN algorithm. However, the $q$-Jacobian matrix tends to be inaccurate when the initial point is set far from the solution since its approximation is based on the assumption that the SSE is small. Otherwise, when $\lambda(k) \to \infty$, the term $\left(\mathbf{J}_q^\mathrm{T}\mathbf{J}_q + \lambda^{(k)}\mathbf{I}\right)^{-1} \approx \mathbf{I}$, the $q$-Hessian matrix approximated with $\mathbf{J}_q^\mathrm{T}\mathbf{J}_q$ loses its importance and the $q$-LM algorithm acts as the $q$-SD algorithm.

## 4 Conclusion

In this paper a new LM model based on *q*-calculus is presented. The *q*-LM algorithm is an iterative procedure which uses damping factor $\lambda$ to alternate between the *q*-SD and the *q*-GN algorithms to minimize the distance between current and previous iterate and the SSE of the corresponding values of the objective function. If $q \to 1$ the *q*-LM algorithm tends to the LM algorithm. When $\lambda \to 0$ the *q*-LM algorithm acts like the *q*-GN algorithm which speeds up the convergence near the stationary point. When $\lambda \to \infty$ the algorithm acts as the *q*-SD algorithm which quickly approaches the solution from the distance.

## References

[1] Michael Lampton. Damping-undamping strategies for the Levenberg-Marquardt least-squares method. Computers in phisics, vol. 11, no. 1, pp. 110-115, 1997.

[2] Henri Gavin. The Levenberg-Marquardt algorithm for nonlinear Least Squares Curve Fitting Problems. Duke University: Department of Civil and Environmental Engineering. Spetember 18, 2020.

[3] Aline Cristina Sotteroni, Roberto Luiz Galski and Fernando Manuel Ramos. The q-gradient method for continuous global optimization. In AIP Conference Proceedings 1558, pp.2389-2393, 2013.

[4] Alfonso Croeze, Lindsey Pittman and Winnie Reynolds. Solving nonlinear least squares problems with Gauss-Newton and Levenberg-Marquardt methods. Retrieved April 26, 2021 from https://www.math.lsu.edu_system_files_MunozGroup1-presentation.

[5] Marcel Maeder, Nichola McCann and S. Norman. Model-based data fitting. In Comprehensive Chemometrics: Chemical and Biochemical Data Analysis, volume 3, pp. 413-436 http://dx.doi.org/10.1016/B978-044452701-1.00058-2, 2009.

[6] Aline Cristina Sotteroni, Roberto Luiz Galski and Fernando Manuel Ramos. The q-gradient vector for unconstrained continuous optimization problems. In B. Hu, et al. (eds.). Operation research proceedings 2010, Springer-Verlag Berlin Heidelberg, pp.365- 370, 2011.

[7] Slađana Marinković, Predrag Rajković and Miomir Stanković. The q-iterative methods in numerical solving of some equations with infinite products. In Facta Universitats, Ser. Math. Inform., vol. 28, no. 4, pp.379-392, 2013.

[8] Kin Keung Lai, Shashi Kant Mishra, Geetanjali Panda, Suvra Kanti Chakraborty, Mohammad Esmael Samei and Bhagwat Ram. A limited memory q-BFGS algorithm for unconstrained optimization problems. Journal of Applied Mathematics and Computing, published online 08 September 2020, https://doi.org/10.1007/s12190-020-01432-6, 2020.

[9] Danijela Protić and Miomir Stanković. The q-Gauss-Newton method for unconstrained nonlinear optimization. Available: https://arxiv.org/ftp/arxiv/papers/2105/2105.12994.pdf.

[10] Kin Keung Lai, Shashi Kant Mishra, Geetanjali Panda, Md Abu Talhamainuddin Ansary and Bhagwat Ram. On q-steepest descent method for unconstrained multiobjective optimization problems. In AIMS Mathematics, vol. 5, no. 6, pp.5521-5540, 2020.

[11] Stanford ENGR 108. Intro to Applied Linear Algebra. Lecture 51. Levenberg Marquardt, 18 march 2021.

[12] Victor Kac and Pokman Cheung. Quantum calculus. Springer, New York, 2002.

[13] Shashi Kant Mishra and Bhagwat Ram. Introduction to unconstrained optimization with R. Springer, https://doi.org/10.1007/978-981-15-0894-3_9 , pp.132-246, 2019.

[14] J. Koekoev and Roelof Koekoev. A note on q-derivative operator. In Journal of Mathematical Analysis and Applications, vol 176, pp.627-634, 1993.

[15] Slađana Marinković, Predrag Rajković and Miomir Stanković. On q-Newton-Kantorowich method for solving systems of equations. In: Applied Mathematics and Computation., vol. 168, pp.1432-1448, 2005.

[16] Pierre Courrieu. Fast Computation of Moore-Penrose inverse matrices. Neural Information Processing – Letters and Reviews, vol.8, no. 2, August 2005, pp.25-29.

[17] Ross MacAusland. The Moore-Penrose inverse and least squares. MATH 420: Advanced Topics in Linear Algebra. April 16, 2014, pp.1-10.